\documentclass{amsart}

\newtheorem{theorem}{Theorem}[section]

\theoremstyle{definition}

\theoremstyle{remark}

\newcommand{\R}{{{\mathbb R}}}

\numberwithin{equation}{section}

\begin{document}

\title{A sufficient condition for finite time blow up of the nonlinear Klein-Gordon equations with arbitrarily positive initial energy}

\author{Yanjin Wang}
\address{Graduate School
of Mathematics, University of Tokyo, 3-8-1 Komaba, Meguro, Tokyo,
153-8914, Japan.}
\email{wangyj@ms.u-tokyo.ac.jp}
\thanks{The author wishes to express his deep
gratitude to Prof. Hitoshi Kitada for his constant encouragement and
kind guidance. The study is supported by Japanese Government
Scholarship.}


%    General info


\begin{abstract}
In this paper we consider the nonexistence of global solutions of
a Klein-Gordon equation of the form
\begin{eqnarray*}
u_{tt}-\Delta u+m^2u=f(u)& (t,x)\in [0,T)\times\R^n.
\end{eqnarray*}
Here $m\neq 0$ and the nonlinear power $f(u)$ satisfies some
assumptions which will be stated later. We give a sufficient
condition on the initial datum with arbitrarily high initial energy
such that the solution of the above Klein-Gordon equation blows up
in a finite time.
\end{abstract}

\maketitle

\section{Introduction}
This paper studies the nonexistence of global solutions of the
nonlinear Klein-Gordon equation
\begin{eqnarray}
\left\{\begin{array}{r@{,\quad  }l}
 u_{tt}-\Delta u+m^2u=f(u)& (t,x)\in [0,T)\times \R^n\\
 u(0,x)=u_0(x)& x\in\R^n\\
 u_t(0,x)=u_1(x)&x\in\R^n
\end{array}
\right.\label{KGE}
\end{eqnarray}
where $\Delta$ is Laplacian operator on $\R^n$ ($n\geq 1$),
$u_0(x)$ and $u_1(x)$ are real valued functions, $m\neq 0$ is a real
constant, the right hand side $f(u)$ is a real-valued, nonlinear function of $u$. Without loss
of generality, we may assume $m=1$ throughout the paper.

The above problem has various applications in the area of
nonlinear optics, plasma physics, fluid mechanics, etc. There is a
literature on the Cauchy problem for the equation (\ref{KGE})
(see, for instance,
\cite{CCazenave}\cite{CGinibre1}\cite{CGinibre2}\cite{CKlainerman}\cite{CPeter}\cite{CSimon}\cite{BStrauss}\cite{Zhang}
and the papers cited therein). The results about the blow up
properties for the local solution of the equation (\ref{KGE}) are
investigated by
\cite{BBall}\cite{BLevine1}\cite{BPayne}\cite{BStrauss}\cite{Zhang}.
For the typical form of the nonlinear power,
$f(u)=|u|^{p-1}u$ with $1<p<\frac{n+2}{n-2}$ ($n\geq 2$), we refer
the result \cite{Zhang} to readers, which established a sharp
condition for the global existence and blow up of the solution of
the Klein-Gordon equation (\ref{KGE}). But in \cite{Zhang}, the
initial energy had an upper bound. This motivates us to consider
the problem how about the solution when the initial energy is
arbitrarily large.

In the paper we investigate the above problem for the equation
(\ref{KGE}) with the nonlinear power $f(u)$ which satisfies that there
exists $\epsilon>0$ such that for any $s\in \R$,
\begin{eqnarray}
f(s)s\geq (2+\epsilon) F(s),\label{PC1}
\end{eqnarray}
where $\displaystyle F(s)=\int_0^s f(\xi)d\xi$.

 For the above
nonlinear power, the nonexistence of the global
solution of an abstract wave equation (including the Klein-Gordon
equation) was established in \cite{BLevine1} when the initial energy was negative.
However until recently there has been very little work on the
global existence and global nonexistence of the solutions of the
initial value problem for Klein-Gordon equations when the initial
energy is arbitrarily high.

The purpose of this paper is to establish a sufficient condition
of the initial datum with arbitrarily high initial energy such
that the corresponding local solution of the nonlinear
Klein-Gordon equation (\ref{KGE}) blows up in a finite time. As
far as we know, this is the first blow up result for the
Klein-Gordon equations with arbitrarily high initial energy on the
whole space $\R^n$. Our proof is based on the concavity method
which firstly introduced by Levine \cite{BLevine1}\cite{Levine2}.
We note that our proof is very simple, and it can be applied to
some other equations.

Here we also refer the result \cite{BLevine3} to readers. Levine
and Todorova \cite{BLevine3} studied a damped wave equation in the
following form
\begin{eqnarray}
\left\{
\begin{array}{l}
u_{tt}+a|u_t|^{\rho-1}u_t-\Delta u=b|u|^{p-1}u-q(x)^2u \mathrm{\
in\ }[0, T)\times \R^n\\
u(0,x)=u_0(x) \mathrm{\ in\ } \R^n\\
u_t(0,x)=u_1(x) \mathrm{\ in\ } \R^n
\end{array}
\right.
 \label{WE}
\end{eqnarray}
where $a, b>0$, $\rho\geq 1$, and $q(x)^2\geq 0$. They proved the
existence of the initial datum with arbitrarily high initial
energy such that the solution blows up in a finite time when
$1\leq\rho<p$. More recently, for the case, $\rho=1$ and $q(x)=0$,
when the initial energy is arbitrarily high, Gazzola and Squassina
\cite{Gazzola} have obtained the blow up result on an open bounded
Lipschitz subset of $\R^n$. We note that their proof cannot be
extended to the whole space $\R^n$. Our proof can be easily
adapted to the damped wave equation (\ref{WE}) with $\rho=1$ and
$q(x)=m\neq 0$. In the paper we will not discuss it in detail but
will make some remarks.

The paper is organized as follows. Section 2 introduces
some notations and known result, and states our main result.
In Section 3 we prove the main result based on a concavity
argument. We make some remarks in the last section.

\section{Preliminary and main result}
Before we state the main result, we introduce some
notations. We denote by $\|\cdot\|_q$ the $L^q(\R^n)$
norm for $1\leq q\leq \infty$, and we define the spaces: $H^1(\R^n)=\{u\in
L^2(\R^n);\ \|u\|_{H^1(\R^n)}=\|(1-\Delta)^{1/2}u\|<\infty\}$, and $H_0^1(\R^n)=\{u\in H^1(\R^n);\ \mathrm{supp}(u)$ $\mathrm{is\
compact\ in\ }\R^n\}$. For simplicity we will denote
$\int_{\R^n}$ by $\int$. The notation $t\rightarrow T^-$ means $t< T$ and
$t\rightarrow T$.

We first state the result of the local existence established in
\cite{CCazenave}
\begin{theorem}
Let the initial data $(u_0, u_1)\in H^1_0(\R^n)\times L^2(\R^n)$,
and let $f$ satisfy the following conditions: $f(0)=0$ and
\begin{eqnarray}
|f(\lambda_1)-f(\lambda_2)|\leq
c(|\lambda_1|^{p-1}+|\lambda_2|^{p-1})|\lambda_1-\lambda_2|
\end{eqnarray}
for all $\lambda_1, \lambda_2\in \R$ and some constant $c>0$, and
\begin{eqnarray}
1<p<\frac{n}{n-2} \mathrm{\ if\ } n\geq 3;\quad 1<p<\infty \mathrm{\ if\
} n=1, 2.
\end{eqnarray}
Then there is a unique local solution $u(t,x)$ of the equation
(\ref{KGE}) on a maximal time interval $[0, T_{\max})$ satisfying
$u(0,x)=u_0(x)$ and $u_t(0,x)=u_1(x)$. In addition, $u(t,x)$
satisfies
\begin{eqnarray}
E(t)=E(0)
\end{eqnarray}
where
\begin{eqnarray}
E(t)=\frac{1}{2}\int\left(|u_t(t,x)|^2+|u(t,x)|^2+|\nabla
u(t,x)|^2-2F(u(t,x))\right)dx.\label{CE}
\end{eqnarray}
\end{theorem}

In order to state the main result, we next define a function,
\begin{eqnarray}
I(u)=\int (|u(x)|^2+|\nabla u(x)|^2)dx-\int f(u(x))u(x) dx.
\end{eqnarray}

Now we are in the position to state our main result:
\begin{theorem} Let $f(s)$ satisfy the assumption
(\ref{PC1}).
 If the initial datum $(u_0,u_1)\in H^1_0(\R^n)\times
L^2(\R^n)$ satisfies that
\begin{eqnarray}
&&E(0)>0, \label{TC1}\\
&&\|u_0\|^2\geq \frac{2(2+\epsilon)}{\epsilon}E(0),\label{TC2}\\
&&I(u_0)<0,\label{TC3}\\
&&\int u_0u_1dx> 0,\label{TC4}
\end{eqnarray}
where $\epsilon>0$ is the constant in (\ref{PC1}), then the
corresponding local solution $u(t,x)$ of the equation (\ref{KGE})
will blow up in a finite time $T_{\max}<\infty$, that is to say,
\begin{eqnarray}
\lim_{t\rightarrow T_{\max}^-}\|u(t,\cdot)\|^2\rightarrow\infty.
\end{eqnarray}
\end{theorem}

\section{Proof of Theorem 2.2}
The proof is split into two steps.

We first prove that
\begin{eqnarray}
&&\|u(t,\cdot)\|^2\geq \frac{2(2+\epsilon)}{\epsilon}E(0),\label{(3.1)}\\
&&I(u(t,\cdot))<0\label{(3.2)}
\end{eqnarray}
for every $t\in[0,T_{\max})$.

We prove $I(u(t,\cdot))<0$ for every $t\in[0, T_{\max})$.
Suppose to the contrary that there exists a time $T>0$ such that
\begin{eqnarray}
T=\min\{t\in(0, T_{\max}); I(u(t,\cdot))=0\}.\nonumber
\end{eqnarray}
We now define an auxiliary function:
\begin{eqnarray}
G(t)=\int|u(t,x)|^2dx.\label{AFG}
\end{eqnarray}
By simple computations we have
\begin{eqnarray}
&&G^\prime(t)=2\int uu_tdx,\\
&&\frac{1}{2}G^{\prime\prime}(t)=\int |u_t|^2dx +\int(
f(u)u-|\nabla u|^2-|u|^2)dx.\label{GTT}
\end{eqnarray}
Noting the assumption $I(u(t,\cdot))<0$ for every $t\in[0, T)$, we
obtain $G^{\prime\prime}(t)>0$ for every $t\in[0,T)$, which
implies that $G^\prime(t)$ is strictly increasing on $[0,T)$, then
by (\ref{TC4}) we see $G^\prime(t)>0$ for every $t\in[0,T)$. In
other words, we obtain that $G(t)$ is also strictly increasing on
$[0, T)$. Thus we obtain
\begin{eqnarray}
G(t)>G(0)\geq\frac{2(2+\epsilon)}{\epsilon}E(0)\nonumber
\end{eqnarray}
for every $t\in (0, T)$. From the continuity of $u(t,x)$ at $t=T$ it
follows that
\begin{eqnarray}
G(T)=\|u(T,\cdot)\|^2>\frac{2(2+\epsilon)}{\epsilon}E(0).\label{Contradiction1}
\end{eqnarray}
But, by (\ref{CE}) it is obvious that
\begin{eqnarray}
\|u(T,\cdot)\|^2+\|\nabla u(T,\cdot)\|^2-2\int F(u(T,x))dx\leq
2E(T)=2E(0).\nonumber
\end{eqnarray}
Noting the assumption $I(u(T,\cdot))=0$ and (\ref{PC1}), we have
\begin{eqnarray*}
\|u(T,\cdot)\|^2+\|\nabla u(T,\cdot)\|^2\geq (2+\epsilon)\int
F(u(T,x))dx,
\end{eqnarray*}
 we now
obtain
\begin{eqnarray}
\|u(T,\cdot)\|^2+\|\nabla u(T,\cdot)\|^2\leq
\frac{2(2+\epsilon)}{\epsilon}E(0).\label{Contradiction2}
\end{eqnarray}
Obviously there is a contradiction between (\ref{Contradiction1})
and (\ref{Contradiction2}). Thus we have proved that for every
$t\in[0, T_{\max})$
\begin{eqnarray}
I(u(t,\cdot))<0.\label{ILeq0}
\end{eqnarray}

By the proof above, we also see that, if
$I(u(t,\cdot))<0$ for every $t\in[0, T_{\max})$ then $G(t)$ is
strictly increasing on $[0, T_{\max})$. Thus (\ref{ILeq0}) implies
\begin{eqnarray}
\|u(t,\cdot)\|^2> \frac{2(2+\epsilon)}{\epsilon}E(0)\label{(3.9)}
\end{eqnarray}
for every $t\in(0, T_{\max})$.

 Therefore we have completed the proof of (\ref{(3.1)}) and (\ref{(3.2)}). We next prove the blow up result of the
equation (\ref{KGE}).

By (\ref{GTT}) and (\ref{(3.9)}), we see that
\begin{eqnarray*}
\frac{1}{2}G^{\prime\prime}(t)&=&\int |u_t|^2dx-\int(|u^2|+|\nabla
u|^2)dx+\int f(u)u dx\\
&\geq&(2+\frac{\epsilon}{2})\int|u_t|^2dx+\frac{\epsilon}{2}\int(|u|^2+|\nabla
u|^2)dx-(2+\epsilon)E(0)\\
&\geq& (2+\frac{\epsilon}{2})\int |u_t|^2dx.
\end{eqnarray*}
By the Schwarz inequality, we therefore have
\begin{eqnarray}
\ \\
G^{\prime\prime}(t)G(t)-\frac{4+\epsilon}{4}(G^\prime(t))^2&\geq&(4+\epsilon)\left\{\int|u|^2dx\int|u_t|^2dx-\left(\int
u u_tdx\right)^2\right\}\nonumber\\
&\geq& 0\nonumber
\end{eqnarray}
for every $t\in[0,T_{\max})$. As $\displaystyle\frac{4+\epsilon}{4}>1$, letting
$\displaystyle\alpha=\frac{\epsilon}{4}$, we have
\begin{eqnarray}
(G^{-\alpha})^{\prime}&=&-\alpha G^{-\alpha-1}G^\prime(t)<0,\\
(G^{-\alpha})^{\prime\prime}&=&-\alpha
G^{-\alpha-2}\left[G^{\prime\prime}(t)G(t)-\frac{4+\epsilon}{4}(G^\prime(t))^2\right]\label{GC}\\
&\leq& 0\nonumber
\end{eqnarray}
for every $t\in(0, T_{\max})$, which means that the function
$G^{-\alpha}$ is concave. Obviously $G(0)>0$,
then from (\ref{GC}) it follows that the function
$G^{-\alpha}\rightarrow 0$ when $t<T_{\max}$ and $\displaystyle
t\rightarrow T_{\max}$ ($\displaystyle T_{\max}<\frac{\epsilon
G(0)}{4G^\prime(0)}$). Therefore, we see that there exists a
finite time $T_{\max}>0$ such that
\begin{eqnarray}
\lim_{t\rightarrow T_{\max}^-}\|u(t,\cdot)\|^2\rightarrow\infty.
\end{eqnarray}
Thus, the proof of Theorem 2.2 is completed.

\section{Some remarks}
First we consider a special case of the nonlinear power $f(u)$,
which is the typical form of the nonlinear power,
$f(u)=|u|^{p-1}u$ with $1<p<\frac{n+2}{n-2}$ and $n\geq 2$. We
define the following action
\begin{eqnarray}
J(u)=\frac{1}{2}\int\left(|u|^2+|\nabla
u|^2-\frac{2}{p+1}|u|^{p+1}\right)dx.
\end{eqnarray}
As in \cite{Berestycki}, we have that
\begin{theorem}
There exists $\bar{u}(x)\in M$, such that
\begin{eqnarray}
d=\inf_{u\in M} J(u)=J(\bar{u})
\end{eqnarray}
where $M=\{u\in H^1(\R^n); I(u)=0, u\neq 0\}$. Moreover, $\bar{u}$
is the ground state of the following elliptic equation
\begin{eqnarray*}
-\Delta u+u=u^p, \quad u\in H^1(\R^n).
\end{eqnarray*}
\end{theorem}
Thus, when $E(0)<d$, as in the first part of the proof of Theorem 2.2, we can prove that if there exists $T$ such that $I(u(T,\cdot))=0$, the function $u(T,x)\in H^1_0(\R^n)\backslash \{0\}$ satisfies
\begin{eqnarray}
&&I(u(T,\cdot))=0,\\
&&J(u(T,\cdot))\leq E(T)=E(0)<d.
\end{eqnarray}
This contradicts Theorem 4.1. Thus for the case $E(0)<d$ the blow-up result holds. This reproduces the result in \cite{Zhang}\cite{Ohta}.

Secondly, we consider an example to which our method is
applicable. For the damped wave equation (\ref{WE}) with $\rho=1$
and $q(x)^2=m^2=1$, we define the auxiliary function $G(t)$
corresponding to (\ref{AFG}) in the following form
\begin{eqnarray}
G(t)=\int |u(t,x)|^2dx +\int_0^t\|u(\tau,\cdot)\|^2d\tau +
(T_0-t)\|u_0\|^2,\nonumber
\end{eqnarray}
where $T_0>0$ is some constant. Then the argument in
\cite{Gazzola} proves $I(u(t))<0$. Thus similarly to the proof of
our Theorem 2.2, we get the blow up result for the damped wave
equation.

\bibliographystyle{amsplain}

\end{document}